\documentclass{article}

\usepackage[T2A]{fontenc}
\usepackage[utf8]{inputenc}
\usepackage[russian,english]{babel}
\usepackage[tbtags]{amsmath}
\usepackage{amsfonts,amssymb}
\usepackage{array}
\usepackage{hyperref}
\usepackage{breakurl}
\usepackage{mathrsfs}

\begin{document}

\author{Yu.\,V.\,Matiyasevich}
\title{Yet another representation for the sum of reciprocals
of the nontrivial zeros of the Riemann zeta-function\footnote{The author
is very grateful to  Peter Zvengrowski (University of
Calgary) for his help with the English.}}

\maketitle

The famous Riemann Hypothesis is
a statement about positions of the zeroes of Riemann's zeta function,
 which can be defined for $\Re(s)>1$ by the Dirichlet series
\begin{equation}
  \zeta(s)=\sum_{n=1}^\infty n^{-s}.
\label{zetadir}\end{equation}
Being analytic, this function is also uniquely defined by its
expansion into
  Laurent series
\begin{equation}
  \zeta(s)=\frac{1}{s-1}+\sum_{n=0}^\infty \gamma_n\frac{(s-1)^n}{n!}.
\label{zetalor}\end{equation}
Thus the Riemann Hypothesis is a statement about the
 \emph{infinite} sequence of real numbers $\gamma_0,\gamma_1,\dots$,
 known as  \emph{Stieltjes constants}
($\gamma_0=\gamma=.577215\dots$ is  the  \emph{Euler constant}).

It is well-known that in order to prove the Riemann Hypothesis it would be
sufficient to establish the validity of a suitable infinite sequence
of polynomial inequalities
\begin{equation}
  P_1(\gamma_0,\dots,\gamma_{m_1})>0,\ \dots,\
  P_n(\gamma_0,\dots,\gamma_{m_n})>0,\ \dots,
\label{Pn}\end{equation}
each of which contains only  \emph{finitely} many
Stieltjes constants  (and, possibly, some other classical
constants), hence allowing numerical verification.

There are many ways to select polynomials $P_n$ giving such a reformulation of
the Riemann Hypothesis. In one of fairly well-known ways described in
 \cite{li} (cf. \cite{bomb}),
\begin{equation}
 P_n(\gamma_0,\dots,\gamma_{n-1})=
 \sum_{\rho}\
 \left(1-\left(1-\frac{1}{\rho}\right)^{n}\right),
\label{Li}\end{equation}
where the the summation is taken over all non-trivial zeroes
of the zeta function.  Earlier, the author considered in
 \cite{me1}
another choice of polynomials  $P_n$ for which
\begin{equation}
  P_n(\gamma_0,\dots,\gamma_{n^2-1})=
  \sum_{j_1=1}^\infty\dots\sum_{j_n=1}^\infty G_n(\rho_{j_1},\dots,\rho_{j_n}),
\label{Grom}\end{equation}
where
\begin{equation}
 G_n(z_1,\dots,z_n)=
 \prod_{j=1}^n\frac{1}{z_j(1-z_j)}\prod_{k=j+1}^n{
 \left(\frac{1}{z_j(1-z_j)}-\frac{1}{z_k(1-z_k)}\right)^2},
\label{G}\end{equation}
and $\rho_1,\rho_2,\dots$ is the sequence
of zeta zeroes with positive imaginary parts.

So far nobody was able to prove all the required
inequalities \eqref{Pn} for any choice of po\-ly\-no\-mi\-als~$P_n$.
The difficulty might be explained by the way in which the Stieltjes constants
are defined. The original definition via formal differentiation of
 \eqref{zetadir} gives slowly convergent series:
\begin{equation}
  \gamma_m=\sum_{n=1}^\infty \left(\frac{\ln^m(n)}{n}-
  \int_n^{n+1}\frac{\ln^m(t)}{t}d t\right).
\label{stiel}\end{equation}

Substitution of  \eqref{stiel} into $P_n(\gamma_0,\dots,\gamma_{m_n})$
results in a multisum, and it is difficult to determine its sign.
However, besides  \eqref{stiel}, different authors found a number of
other representations for the Stieltjes constants in the form
of infinite sum or integral (cf., for example  \cite{coffey} and references
there).
The author finds it plausible that for certain
choice of polynomials  $P_n$ and representations for $\gamma_0,\ \gamma_1, \dots$,
one could find for $P_n(\gamma_0,\dots,\gamma_n)$ representations in the form
of (multi)sums with positive summands or integrals with positive
 integrands.

The author was able to obtain such representations for
 $n=0$ in \eqref{Li} and in
\eqref{Grom}. In both cases
\begin{multline}
  2P_0(\gamma_0)=2\sum_{j=1}^\infty\left(\frac{1}{\rho_j}+\frac{1}{1-\rho_j}
  \right)=
 2\sum_{j=1}^\infty\frac{1}{\rho_j(1-\rho_j)}=\\=
 \gamma-\ln(\pi)-2\ln(2)+2=0.0461914\dots
\label{P0}\end{multline}
From \cite{me1}
\begin{equation}
  \gamma-\ln(4\pi)+2=\int_1^\infty\frac{1-\{q\}^2}{4q^2(q+1)^2}\,dq,
\label{P01}\end{equation}
where $\{q\}$ is the fractional part of  $q$, so the
positivity of the integrand is evident.
From
\cite{me2}
\begin{equation}
  \gamma-\ln(4\pi)+2=\sum_{n=1}^\infty\left(
  \psi(n)-\int_{n-1/2}^{n+1/2}\psi(q)\,d q\right),
\label{P12}\end{equation}
where $\psi(q)=\Gamma'(q)/\Gamma(q)$ is the logarithmic derivative of the
gamma-function, called the digamma function. The integrand in \eqref{P12} is
positive because $\psi''''(q)<0$ for $q>1/2$.

The aim of the present note is to give yet another
representation for \eqref{P0} which makes evident its positivity,
namely, the equality \eqref{itog}.

At first sight  \eqref{P0} contains three numbers of different nature --
$\gamma$, $\ln(\pi)$, and $\ln(2)$, and it is not clear what is the ``reason''
for the positivity of $P_0(\gamma_0)$. In order to get a desired representation with
positive summands or  integrands
one needs to find representations for these constants that
would look similar.
The new representation is based on duality, indicated in
\cite{sondow0}, between  $\gamma$ and $\ln(\pi)$:
\begin{equation}
  \gamma=\sum_{n=1}^\infty \left(\frac{1}{n}-\int_n^{n+1}\frac{1}{t}d t\right),
  \qquad
  \ln\left(\frac{4}{\pi}\right)=\sum_{n=1}^\infty(-1)^{n-1}
  \left(\frac{1}{n}-\int_n^{n+1}\frac{1}{t}d t\right)
\label{gammapi}\end{equation}
(the left equality in  \eqref{gammapi} is just the case  $m=0$ of \eqref{stiel}).

In \cite{sondow1} the two equalities from \eqref{gammapi}
were transformed respectively into
\begin{equation}
  \gamma=\sum_{n=2}^\infty (-1)^{n}\frac{
  N_1\left(\lfloor\frac{n}{2}\rfloor\right)+
        N_0\left(\lfloor\frac{n}{2}\rfloor\right)}{n},
    \quad
    \ln\left(\frac{4}{\pi}\right)=
  \sum_{n=2}^\infty (-1)^{n}\frac{
  N_1\left(\lfloor\frac{n}{2}\rfloor\right)-
        N_0\left(\lfloor\frac{n}{2}\rfloor\right)}{n},
\label{vaccagammapi}\end{equation}
where $N_0(m)$ and  $N_1(m)$ denote the number of  zeroes and
units the binary expansion of~$m$. Clearly,
  $N_1\left(\lfloor\frac{n}{2}\rfloor\right)+
        N_0\left(\lfloor\frac{n}{2}\rfloor\right)=
        \lfloor\log_2(n)\rfloor$,
so the left hand side in  \eqref{vaccagammapi}
is just another transcription of
 the so-called \emph{Vacca series} (\cite{vacca}, cf. also \cite{hardy}).
Pairwise grouping of the summands has given dual equalities
\begin{equation}
  \gamma=\sum_{n=1}^\infty \frac{
  N_1\left({n}\right)+
        N_0\left({n}\right)}{2n(2n+1)},
    \qquad
    \ln\left(\frac{4}{\pi}\right)=
  \sum_{n=1}^\infty \frac{
  N_1\left({n}\right)-
        N_0\left({n}\right)}{2n(2n+1)}.
\label{vaccagammapifast}\end{equation}

An infinite series of representations for $\gamma$
was then constructed in \cite{sondow1} by accelerating the convergence of
 the left-hand side equation from
 \eqref{vaccagammapifast}; a particular case, the equality
\begin{equation}
  \gamma=\frac{1}{2}+\sum_{n=1}^\infty \frac{
  N_1\left({n}\right)+
        N_0\left({n}\right)}{2n(2n+1)(2n+2)},
    \label{addison}\end{equation}
was established earlier in  \cite{addison} by a different method.
In a similar way the right-hand side of  \eqref{vaccagammapifast}
 can be transformed into
a representation that is dual to  \eqref{addison}:
\begin{equation}
  \ln\left(\frac{2}{\pi}\right)=-\frac{1}{2}+\sum_{n=1}^\infty \frac{
  N_1\left({n}\right)-
        N_0\left({n}\right)}{2n(2n+1)(2n+2)},
    \label{vaccapifast}\end{equation}
which is essentially  the case
 $B=2$ of Corollary  19 from \cite{pile}.

Summing up \eqref{addison} and \eqref{vaccapifast}, we get the
equality
\begin{equation}
  \gamma-\ln(\pi)+\ln(2)=\sum_{n=1}^\infty \frac{
  2N_1\left({n}\right)}{2n(2n+1)(2n+2)},
    \label{pochti}
\end{equation}
which has all positive summands
in the right hand side..

The value of \eqref{pochti} differs from the desired value of
\eqref{P0} by $3\ln(2)-2$, so we need to find a suitable expression
for   $\ln(2)$. An easy calculation give the equality
 \begin{equation}
  \frac{3}{4}-\ln(2)=\sum_{n=1}^\infty \frac{
  1}{2n(2n+1)(2n+2)},
    \label{log2}\end{equation}
and respectively
\begin{equation}
  \gamma-\ln(\pi)-2\ln(2)+\frac{9}{4}=\sum_{n=1}^\infty \frac{
  2N_1\left({n}\right)+3}{2n(2n+1)(2n+2)}.
    \label{pochtipochti}
\end{equation}
Dropping the first two summands in the right-hand side
we get the desired equality
 \begin{equation}2\sum_{j=1}^\infty\left(\frac{1}{\rho_j}+\frac{1}{1-\rho_j}
  \right)=
  \gamma-\ln(4\pi)+2=\sum_{n=3}^\infty \frac{
  2N_1\left({n}\right)+3}{2n(2n+1)(2n+2)}.
    \label{itog}
\end{equation}

It remains intriguing whether one could get similar representations for\linebreak
 $P_2(\gamma_0,\dots,\gamma_{m_2})$, $P_3(\gamma_0,\dots,\gamma_{m_3})$, \dots
 in the form of (multi)sums with positive summands
 for  $P_n$ defined by
\eqref{Li}, \eqref{Grom}, or in any other way implying
the validity of the Riemann Hypothesis.


\begin{thebibliography}{99}

\bibitem{addison}
A. W. Addison,
A series representation for Euler's constant,
\emph{Amer. Math. Monthly} 74 (1967), 823--824.


\bibitem{bomb}
 E. Bombieri,  J. C. Lagarias, Complements to Li's criterion for
 the Riemann hypothesis,
\emph{J. Number Theory} 77 (1999), 274--287.


\bibitem{coffey}
M. W. Coffey,
Addison-type series representation for the Stieltjes
constants,  \url{http://arxiv.org/abs/0912.2391}, (2009).


\bibitem{hardy}
G. H. Hardy, Note on Dr.\,Vacca's Series for gamma,
\emph{Quart. J. Pure Appl. Math.} 43 ( 1912), 215--216.


\bibitem{li}
X.-J. Li,  The positivity of a sequence of numbers and the Riemann hypothesis,
\emph{J. Number Theory} 65 (1997), 325--333.

\bibitem{me1}
Yu. V. Matiyasevich,
An analytic representation for the sum of reciprocals
of the nontrivial zeros of the Riemann zeta-function
\emph{Proceedings of the Steklov Institute of Mathematics},
 163 (1985),  211--213 (translated from
Ю. В. Матиясевич,
Одно аналитическое представление для суммы величин,
обратных к нетривиальным нулям дзета-функции Римана,
\emph{Труды МИАН}, 163 (1984), 181--182,
\url{http://www.mathnet.ru/links/c9ca7c722c53bbebe971fa4d79daf173/tm2325.pdf}).



\bibitem{me2}
Yu. V. Matiyasevich,
A relationship between certain sums over trivial
and nontrivial zeros of the Riemann zeta-function,
\emph{Mathematical Notes},
 45:(1--2) (1989),  131--135,
 doi: 10.1007/BF01158058 (translated from
Ю. В. Матиясевич,
Связи между некоторыми суммами по тривиальным и
нетривиальным нулям дзета-функции Римана,
\emph{Математические заметки}
45:2 (1989),  65--70).


\bibitem{pile}
K. H. Pilehrood,   T. H. Pilehrood,
Vacca-Type series for values of the
generalized Euler constant function
and its derivative. \emph{Journal of Integer Sequences}, Vol. 13 (2010), article 10.7.3;
\url{https://cs.uwaterloo.ca/journals/JIS/VOL13/Pilehrood/pilehrood2.pdf}


 \bibitem{sondow0}
 J. Sondow, Double integrals for Euler's constant and
 $\ln(\frac{4}{\pi} )$ and an analog of
Hadjicostas's formula, \emph{Amer. Math. Monthly} 112 (2005), 61--65.


\bibitem{sondow1}
J. Sondow,
New Vacca-type rational series for Euler's constant $\gamma$ and
its ``alternating'' analog $\ln\frac4\pi$.
\emph{Additive number theory}, Springer, New York, 2010, 331–340;
\url{http://arxiv.org/abs/math/0508042}.


\bibitem{vacca}
 G. Vacca, A new series for the Eulerian constant $\gamma = 0.577...$,
 \emph{Quart. J. Pure Appl. Math.} 41 (1909--1910), 363--366.


\end{thebibliography}
\end{document}